\documentclass[aos,MSNbibl,dvips]{arximspdf}

\doi{10.1214/11-AOS928}
\volume{39}
\issue{5}
\pubyear{2011}
\firstpage{2243}
\lastpage{2243}

\begin{document}
\begin{frontmatter}

\title{Introduction to the Lehmann special section}
\runtitle{Introduction to the Lehmann special section}

\begin{aug}
\author[A]{\fnms{Peter} \snm{B\"uhlmann}}
\and
\author[B]{\fnms{Tony} \snm{Cai}\corref{}\ead[label=e1]{tcai@wharton.upenn.edu}}
\runauthor{T. Cai and P. Bühlmann}
\address[A]{Seminar for Statistics\\
ETH Zentrum, HG G 10.3\\
R\"{a}mistrasse 101\\
8092 Zurich\\
Switzerland}
\address[B]{Department of Statistics\\
The Wharton School\\
University of Pennsylvania\\
Philadelphia, Pennsylvania 19104\\
USA} 
\end{aug}

\received{\smonth{9} \syear{2011}}



\end{frontmatter}

Erich L. Lehmann died at the age of 91 on September 12, 2009. The statistics community lost a great scientist who made fundamental contributions to many areas of statistics. The Council of the Institute of
Mathematical Statistics has decided that this issue of The Annals of
Statistics should be dedicated to the memory of Erich L. Lehmann.

Erich L. Lehmann was a pioneer in establishing a mathematical
framework and theory for the foundations of modern statistics. His
books \textit{Testing Statistical Hypotheses} (1959) and \textit{Theory of Point
Estimation} (1983) have helped to educate several generations of statisticians
around the globe.
Many avenues of research began with the seminal work of
Erich; these have been further developed and brought to a broad range of
fields in statistics, often by his many former Ph.D. students.

The current Special Issue of \textit{The Annals of Statistics} contains three
invited articles. Javier Rojo  discusses Erich's scientific
achievements and provides complete lists of his scientific writings and his
former Ph.D. students. Willem van Zwet  describes aspects of Erich's
life and work, enriched with personal and interesting anecdotes of Erich's long and productive scientific journey. Finally, Peter Bickel, Aiyou Chen and Elizaveta
Levina present a research paper on network models: they dedicate their
contribution to Erich, emphasizing that their new nonparametric
method and issues about optimality have been very much influenced by
Erich's thinking.

Erich Lehmann was not only an outstanding scientist but also a very kind and generous
person. Although he is not among us any more, his scientific legacy
and wonderful memories will always remain.

\printaddresses

\end{document}